# ASYMPTOTIC LOCAL EFFICIENCY OF CRAMÉR–VON MISES TESTS FOR MULTIVARIATE INDEPENDENCE

BY CHRISTIAN GENEST[1], JEAN-FRANÇOIS QUESSY AND
BRUNO RÉMILLARD[1]

*Université Laval, Université du Québec à Trois-Rivières and HEC Montréal*

Deheuvels [*J. Multivariate Anal.* **11** (1981) 102–113] and Genest and Rémillard [*Test* **13** (2004) 335–369] have shown that powerful rank tests of multivariate independence can be based on combinations of asymptotically independent Cramér–von Mises statistics derived from a Möbius decomposition of the empirical copula process. A result on the large-sample behavior of this process under contiguous sequences of alternatives is used here to give a representation of the limiting distribution of such test statistics and to compute their relative local asymptotic efficiency. Local power curves and asymptotic relative efficiencies are compared under familiar classes of copula alternatives.

**1. Introduction.** In a seminal paper concerned with testing the null hypothesis of independence between the $d \geq 2$ components of a multivariate vector with continuous distribution $H$ and marginals $F_1, \ldots, F_d$, Blum, Kiefer and Rosenblatt [1] investigated the use of a Cramér–von Mises statistic derived from the process

$$\mathbb{H}_n(x) = \sqrt{n}\left\{H_n(x) - \prod_{j=1}^{d} F_{j,n}(x_j)\right\}, \qquad x = (x_1, \ldots, x_d) \in \mathbb{R}^d,$$

that measures the difference between the empirical distribution function $H_n$ of $H$ and the product of the marginal empirical distributions $F_{j,n}$ associated with the components of the random vector. As Hoeffding [17] had

Received February 2005; revised March 2006.
[1]Supported by grants from the Natural Sciences and Engineering Research Council of Canada and from the Fonds québécois de la recherche sur la nature et les technologies.

already noted, the asymptotic distribution of this test statistic is generally not tractable, and hence tables of critical values are required for its use. Such tables were provided by Blum, Kiefer and Rosenblatt [1] themselves in the case $d = 2$. They were later expanded to $d \geq 3$ by Csörgő [4] and by Cotterill and Csörgő [2, 3], based on strong approximations of $\mathbb{H}_n$. A bootstrap approach has also been proposed by Jing and Zhu [18].

Despite the anticipation that the Cramér–von Mises statistic $\int \mathbb{H}_n^2 \, dH_n$ should be powerful, most subsequent research focused on the case $d = 2$, where alternative tests, typically based on moment characterizations of independence, were proposed, for example, by Feuerverger [9], Shih and Louis [21], Gieser and Randles [15] and Kallenberg and Ledwina [19].

Curiously, the literature seems to have largely ignored a suggestion of Blum, Kiefer and Rosenblatt [1] to circumvent the inconvenience caused by the complex nature of the limiting distribution of $\mathbb{H}_n$. To be specific, let $X_1 = (X_{11}, \ldots, X_{1d}), \ldots, X_n = (X_{n1}, \ldots, X_{nd})$ be a random sample from the distribution $H$ and for arbitrary $A \subset \mathcal{S}_d = \{1, \ldots, d\}$ with $|A| > 1$, consider

$$G_{A,n}(x) = \frac{1}{\sqrt{n}} \sum_{i=1}^n \prod_{j \in A} \{\mathbf{1}(X_{ij} \leq x_j) - F_{j,n}(x_j)\}.$$

Using the multinomial identity

$$\prod_{j=1}^d (a_j + b_j) = \sum_{A \subset \mathcal{S}_d} \prod_{j \in A} a_j \prod_{j \notin A} b_j, \qquad a_1, \ldots, a_d, \ b_1, \ldots, b_d \in \mathbb{R},$$

with the convention that a product over an empty set of terms equals 1, Blum, Kiefer and Rosenblatt [1] showed that $\mathbb{H}_n$ may be conveniently expressed as

$$\mathbb{H}_n(x) = \sum_{A \subset \mathcal{S}_d, |A| > 1} G_{A,n}(x) \prod_{j \in \mathcal{S}_d \setminus A} F_{j,n}(x_j).$$

Although their paper focused on the case $d = 3$, these authors claimed, as later confirmed by Dugué [8], that under independence, $G_{A,n}$ converges weakly to a continuous centered Gaussian process with covariance function

$$\mathrm{cov}_A(x, x') = \prod_{j \in A} [\min\{F_j(x_j), F_j(x'_j)\} - F_j(x_j) F_j(x'_j)]$$

and whose eigenvalues, given by $\pi^{-2|A|} (i_1 \cdots i_{|A|})^{-2}$ for $i_1, \ldots, i_{|A|} \in \mathbb{N}$, may be deduced from the Karhunen–Loève decomposition of the Brownian bridge. More importantly still, Blum, Kiefer and Rosenblatt [1] and Dugué [8] pointed out that the processes $G_{A,n}$ and $G_{A',n}$ are asymptotically independent whenever $A \neq A'$ so that suitable combinations of statistics based on the individual $G_{A,n}$ processes could be used to test independence.

An obvious limitation of tests based on this approach, however, is the dependence of the asymptotic null distribution of the $G_{A,n}$ on the marginals of $H$. To alleviate this problem, Deheuvels [6] suggested that the original observations $X_1, \ldots, X_n$ be replaced by their associated rank vectors $R_1 = (R_{11}, \ldots, R_{1d}), \ldots, R_n = (R_{n1}, \ldots, R_{nd})$, where

$$R_{ij} = \sum_{\ell=1}^{n} \mathbf{1}(X_{\ell j} \leq X_{ij}), \qquad 1 \leq i \leq n, \ 1 \leq j \leq d.$$

Deheuvels then went on to characterize the asymptotic null behavior of a Möbius decomposition of the copula process

$$(1.1) \qquad \mathbb{C}_n(u) = \sqrt{n}\left\{C_n(u) - \prod_{j=1}^{d} u_j\right\},$$

where the empirical copula [5], defined by

$$(1.2) \qquad C_n(u_1, \ldots, u_d) = \frac{1}{n}\sum_{i=1}^{n}\prod_{j=1}^{d}\mathbf{1}(R_{ij} \leq nu_j),$$

is an estimate of the unique copula $C$ defined implicitly by $C\{F_1(x_1), \ldots, F_d(x_d)\} = H(x_1, \ldots, x_d)$ for all $x_1, \ldots, x_d \in \mathbb{R}$. As the latter reduces to $C(u_1, \ldots, u_d) = u_1 \cdots u_d$ under independence, Deheuvels [6] proposed that this hypothesis be tested using a Cramér–von Mises statistic based on a decomposition of $\mathbb{C}_n$, but he concentrated his own efforts on the identification of the asymptotic null distribution.

Recently, Genest and Rémillard [13] showed how to compute the quantiles for the finite-sample and asymptotic null distribution of Cramér–von Mises statistics based on $\mathbb{C}_n$. Furthermore, they investigated how the $2^d - d - 1$ statistics derived from the rank analogues $\mathbb{G}_{A,n}$ of the $G_{A,n}$ could be combined to obtain a global statistic for testing independence.

This paper enhances the work of Genest and Rémillard [13] by comparing the power of Cramér–von Mises tests of independence based on the copula process $\mathbb{C}_n$ and of those based on four different combination recipes for the terms of its Möbius decomposition. To this end, the local asymptotic behavior of the copula process $\mathbb{C}_n$ is characterized in Section 2 under sequences of alternatives contiguous to independence. Examples of such sequences are considered in Section 3, where it is shown that Clayton's gamma frailty model shares with the equicorrelated Gaussian model the surprising property that the limiting behavior of $\mathbb{G}_{A,n}$ is independent of the sequence of contiguous alternatives for all $A \in \mathcal{S}_d$ with $|A| > 2$. In Section 4, some test statistics for the null hypothesis of independence are considered and their asymptotic behavior under contiguous alternatives is specified. Local power functions are then computed in Section 5. Finally, Section 6 gives additional comparisons between the test statistics using an extension of Pitman's local asymptotic efficiency considered by Genest, Quessy and Rémillard [12] in a bivariate setting.

**2. Limiting behavior of $\mathbb{C}_n$ under contiguous alternatives.** Let $\Theta \subset \mathbb{R}$ be a closed interval and $\mathcal{C} = \{C_\theta : \theta \in \Theta\}$ be a given family of copulas that are monotone in $\theta$ with respect to the concordance ordering and for which $\theta_0 \in \Theta$ corresponds to independence. Take $\delta_n \to \delta \in \mathbb{R}$ such that $\theta_n = \theta_0 + \delta_n/\sqrt{n} \in \Theta$ for $n$ sufficiently large.

Let $Q_n$ be the joint distribution function of the random sample $(X_{11}^{(n)}, \ldots, X_{1d}^{(n)}), \ldots, (X_{n1}^{(n)}, \ldots, X_{nd}^{(n)})$ from the distribution function $C_{\theta_n}\{F_1(x_1), \ldots, F_d(x_d)\}$ and let $C_n$ be the empirical copula of these observations computed using formula (1.2). Finally, let $P_n$ be the joint distribution of a sample of the same size under the independence distribution $F_1 \times \cdots \times F_d$.

The limiting behavior of the sequence $(\mathbb{C}_n)$ of empirical copula processes will be determined under the following conditions:

(i) $C_\theta$ is absolutely continuous and its density $c_\theta$ has a square integrable right derivative $\dot{c}_\theta(u) = \partial c_\theta(u)/\partial \theta$ at $\theta = \theta_0$ for each $u = (u_1, \ldots, u_d) \in (0,1)^d$, the latter satisfying

$$\lim_{n \to \infty} \int_{(0,1)^d} \left[ \sqrt{n}\{\sqrt{c_{\theta_n}(u)} - 1\} - \frac{\delta}{2}\dot{c}_{\theta_0}(u) \right]^2 du_d \cdots du_1 = 0;$$

(ii) for every $u = (u_1, \ldots, u_d) \in (0,1)^d$, one has

$$\dot{C}_{\theta_0}(u) = \lim_{\theta \to \theta_0} \frac{\partial}{\partial \theta} C_\theta(u) = \int_0^{u_1} \cdots \int_0^{u_d} \dot{c}_{\theta_0}(v)\, dv_d \cdots dv_1.$$

The following result extends Proposition 1 of Genest, Quessy and Rémillard [12]:

PROPOSITION 2.1. *Suppose that the underlying copula of a given population belongs to a family $\mathcal{C}$ whose members satisfy assumptions* (i) *and* (ii) *above. Then under $Q_n$, the empirical process $\mathbb{C}_n = \sqrt{n}(C_n - C_{\theta_0})$ converges in law in $\mathcal{D}([0,1]^d)$ to a continuous Gaussian process $\mathbb{C} + \delta \dot{C}_{\theta_0}$ with covariance structure $\Lambda$, where for each $u, u' \in [0,1]^d$, $\Lambda(u, u')$ is given by*

$$C_{\theta_0}(u \wedge u') + (d-1)C_{\theta_0}(u)C_{\theta_0}(u') - C_{\theta_0}(u)C_{\theta_0}(u') \sum_{j=1}^d \left( \frac{u_j \wedge u'_j}{u_j u'_j} \right).$$

PROOF. Introduce $U_{ij}^{(n)} = F_j(X_{ij}^{(n)})$ and define

$$\mathbb{A}_n(u) = \frac{1}{\sqrt{n}} \sum_{i=1}^n \left\{ \prod_{j=1}^d \mathbf{1}(U_{ij}^{(n)} \leq u_j) - \prod_{j=1}^d u_j \right\}.$$

Also, let $\Psi_{j,n}(u) = \sum_{i=1}^{n} \mathbf{1}(U_{ij}^{(n)} \leq u)/n$. Since $R_{ij}^{(n)} = n\Psi_{j,n}(U_{ij}^{(n)})$, it follows from equations (1.1) and (1.2) that

$$(2.1) \quad \mathbb{C}_n(u) = \mathbb{A}_n\{\Psi_{1,n}^{-1}(u_1),\ldots,\Psi_{d,n}^{-1}(u_d)\} + \sqrt{n}\left\{\prod_{j=1}^{d}\Psi_{j,n}^{-1}(u_j) - \prod_{j=1}^{d}u_j\right\}.$$

Under assumption (i), an application of Theorem 3.10.12 of van der Vaart and Wellner [23] implies that under $Q_n$, the sequence $(\mathbb{A}_n)$ of processes converges in $\mathcal{D}([0,1]^d)$ to a continuous Gaussian limit of the form $\mathbb{A} + \delta\dot{C}_{\theta_0}$, where $\dot{C}_{\theta_0}$ is defined as in assumption (ii).

As a consequence of this result, one has that under $Q_n$, the univariate process $\mathbb{A}_n(\mathbf{1}, u_j, \mathbf{1}) = \sqrt{n}\{\Psi_{j,n}(u_j) - u_j\}$ converges in $\mathcal{D}([0,1])$ to $\mathbb{A}(\mathbf{1}, u_j, \mathbf{1}) + \delta\dot{C}_{\theta_0}(\mathbf{1}, u_j, \mathbf{1}) = \mathbb{A}(\mathbf{1}, u_j, \mathbf{1})$ since

$$\dot{C}_{\theta_0}(\mathbf{1}, u_j, \mathbf{1}) = \lim_{\theta \to \theta_0} \frac{C_\theta(\mathbf{1}, u_j, \mathbf{1}) - C_{\theta_0}(\mathbf{1}, u_j, \mathbf{1})}{\theta - \theta_0} = \lim_{\theta \to \theta_0} \frac{u_j - u_j}{\theta - \theta_0} = 0.$$

Using identities (11) and (12) from Chapter 3 in Shorack and Wellner [22], one can see that

$$\sup_{u_j \in [0,1]} |\Psi_{j,n}(u_j) - u_j| = \sup_{u_j \in [0,1]} |\Psi_{j,n}^{-1}(u_j) - u_j| \to 0$$

in probability for each $j \in \mathcal{S}_d$, whence $\sqrt{n}\{\Psi_{j,n}^{-1}(u_j) - u_j\} \to -\mathbb{A}(\mathbf{1}, u_j, \mathbf{1})$ in $\mathcal{D}([0,1])$. Finally, the second summand in (2.1) can be rewritten as

$$\sum_{k=1}^{d}\left(\prod_{i=1}^{k-1}u_i\right)\left\{\prod_{j=k+1}^{d}\Psi_{j,n}^{-1}(u_j)\right\}\sqrt{n}\{\Psi_{k,n}^{-1}(u_k) - u_k\}.$$

Accordingly, one may conclude that under $Q_n$, the process $\mathbb{C}_n$ converges in $\mathcal{D}([0,1]^d)$ to $\mathbb{C} + \delta\dot{C}_{\theta_0}$, where

$$\mathbb{C}(u) = \mathbb{A}(u) - \sum_{k=1}^{d}\mathbb{A}(\mathbf{1}, u_k, \mathbf{1})\prod_{j \neq k}u_j$$

is the limiting process of $\mathbb{C}_n$ under $P_n$ that was identified by Gänssler and Stute [10]. Finally, a straightforward computation shows that the limiting covariance function of $\mathbb{C}_n$ is $\Lambda$, given that $\delta\dot{C}_{\theta_0}$ is a deterministic term. $\square$

Deheuvels [6] proposed the decomposition of $\mathbb{C}_n$ into a collection of asymptotically independent, centered Gaussian processes having a simple covariance function under the null hypothesis of independence. Specifically, for each $A \subset \mathcal{S}_d$, define the linear operator $\mathcal{M}_A$ such that

$$\mathbb{G}_{A,n}(u) = \mathcal{M}_A\{\mathbb{C}_n(u)\} = \frac{1}{\sqrt{n}}\sum_{i=1}^{n}\prod_{j \in A}\{\mathbf{1}(R_{ij} \leq nu_j) - u_j\}$$
$$= \sum_{B \subset A}(-1)^{|A \setminus B|}\mathbb{C}_n(u^B)C_{\theta_0}(u^{A \setminus B}),$$

where in general, $u^B = (u_1^B, \ldots, u_d^B)$ with $u_j^B$ equal to $u_j$ when $j \in B$ and equal to 1 otherwise. The result is that under $P_n$, one obtains a decomposition

$$\mathbb{C}_n(u) = \sum_{A \subset \mathcal{S}_d, |A| > 1} \mathbb{G}_{A,n}(u) \prod_{j \in \mathcal{S}_d \setminus A} u_j$$

of the empirical process $\mathbb{C}_n$ into $2^d - d - 1$ subprocesses that converge jointly to a vector of continuous, centered Gaussian processes $\mathbb{G}_A = \mathcal{M}_A(\mathbb{C})$. Furthermore, their asymptotic covariance structure is given by

$$\Gamma_{A,A'}(u, u') = \text{cov}\{\mathbb{G}_A(u), \mathbb{G}_{A'}(u')\} = \mathbf{1}(A = A') \prod_{j \in A} \gamma(u_j, u'_j),$$

where $\gamma(u, u') = u \wedge u' - uu'$. In other words, the $\mathbb{G}_{A,n}$ are asymptotically independent and their limiting covariance is the same as that of a product of independent Brownian bridges.

The next result, which is a straightforward consequence of Proposition 2.1, gives the asymptotic representation of $\mathbb{G}_{A,n}$ under $Q_n$.

COROLLARY 2.2. *Let $\mathcal{C}$ be a given family of copulas whose members satisfy assumptions* (i) *and* (ii) *above. Then under $Q_n$, the empirical processes $\mathbb{G}_{A,n}$ converge jointly in $\mathcal{D}([0,1]^d)$ to a vector of continuous Gaussian processes $\mathbb{G}_A + \delta\mu_A$ with covariance structure $\Gamma_{A,A'}$, where $\mu_A = \mathcal{M}_A(\dot{C}_{\theta_0})$.*

**3. Examples.** The drift term $\mu_A(u)$ identified in Corollary 2.2 can be computed explicitly in many families of copulas. The following result, known as the Möbius inversion formula, will be useful for that purpose.

LEMMA 3.1. *Let $f$ be a function defined on the subsets $B$ of $\mathcal{S}_d$. For any $A \subset \mathcal{S}_d$, set $F(A) = \sum_{B \subset A} f(B)$. Then $f(A) = \sum_{B \subset A} (-1)^{|A \setminus B|} F(B)$.*

3.1. *The multivariate equicorrelated normal copula.* Let $\Phi$ denote the cumulative distribution function of a standard normal random variable. The multivariate normal copula with $d \times d$ correlation matrix $\Sigma = (\sigma_{jk})$ is defined by $C_\Sigma(u) = H_\Sigma\{\Phi^{-1}(u_1), \ldots, \Phi^{-1}(u_d)\}$, where

$$H_\Sigma(x) = \int_{-\infty}^{x_1} \cdots \int_{-\infty}^{x_d} \frac{1}{(2\pi)^{d/2}|\Sigma|^{1/2}} \exp(-y^\top \Sigma^{-1} y/2) \, dy_d \cdots dy_1.$$

Consider the equicorrelated case in which $\sigma_{jj} = 1$ for all $1 \leq j \leq d$ and $\sigma_{jk} = \rho$ for all $j \neq k$. Write $H_\rho = H_\Sigma$ and $\dot{H}_\rho = dH_\rho/d\rho$. Then

$$\left.\frac{d}{d\rho}|\Sigma|^{-1/2}\right|_{\rho=0} = 0 \quad \text{and} \quad \left.-\frac{1}{2}\frac{d}{d\rho}x^\top \Sigma^{-1} x\right|_{\rho=0} = \sum_{j<k} x_j x_k,$$

where the latter identity follows from the fact that

$$x^\top \Sigma^{-1} x = \frac{(1-\rho)^{d-2}}{|\Sigma|}\left[\{(d-2)\rho+1\}\sum_{\ell=1}^d x_\ell^2 - 2\rho \sum_{j<k} x_j x_k\right].$$

Now if $\varphi(t) = d\Phi(t)/dt$, one gets

$$\lim_{\rho \to 0} \dot{H}_\rho(x) = \int_{-\infty}^{x_1} \cdots \int_{-\infty}^{x_d} \left(\sum_{j<k} y_j y_k\right)\left\{\prod_{\ell=1}^d \varphi(y_\ell)\right\} dy_d \cdots dy_1$$

$$= \sum_{j<k}\left\{\prod_{\ell \neq j,k} \Phi(x_\ell)\right\}\varphi(x_j)\varphi(x_k) = H_0(x) \sum_{j<k} \frac{\varphi(x_j)}{\Phi(x_j)}\frac{\varphi(x_k)}{\Phi(x_k)}.$$

Thus,

$$\dot{C}_{\theta_0}(u) = \lim_{\rho \to 0} \dot{H}_\rho\{\Phi^{-1}(u_1),\ldots,\Phi^{-1}(u_d)\}$$

$$= C_{\theta_0}(u) \sum_{j<k} \frac{1}{u_j u_k} \varphi\{\Phi^{-1}(u_j)\}\varphi\{\Phi^{-1}(u_k)\}.$$

An application of the Möbius inversion formula then yields

$$\mu_A(u) = \mathbf{1}(|A|=2) \prod_{j \in A} \varphi\{\Phi^{-1}(u_j)\}.$$

Hence, when looking at the Möbius decomposition of the multivariate normal model, only the $\mathbb{G}_{A,n}$ with $|A|=2$ have a limiting distribution that differs under the null hypothesis and contiguous alternatives. Accordingly, tests of independence should be based only on the latter because the inclusion of functions of $\mathbb{G}_{A,n}$ for any $|A| > 2$ would contribute nothing to the overall power of the procedure. This observation does not come as a total surprise, given that the multivariate Gaussian dependence structure is completely characterized by the pairwise interactions among the variables.

3.2. *One-parameter multivariate Farlie–Gumbel–Morgenstern copula.* A $d$-variate version of this system of copulas is defined for $\theta \in [-1,1]$ by

$$C_\theta(u) = C_{\theta_0}(u) + \theta \prod_{j=1}^d u_j(1-u_j), \qquad u \in [0,1]^d,$$

with $\theta_0 = 0$ corresponding to independence. It follows easily that $\dot{C}_{\theta_0}(u) = u_1(1-u_1) \times \cdots \times u_d(1-u_d)$. Since $\dot{C}_{\theta_0}(u^B) = 0$ whenever $B \neq \mathcal{S}_d$, one gets

$$\mu_A(u) = (-1)^d \dot{C}_{\theta_0}(u)\mathbf{1}(A = \mathcal{S}_d).$$

Thus, in this case, $\mu_A$ vanishes unless $A = \mathcal{S}_d$, implying that in contrast to the multivariate normal case, tests of independence based on the $\mathbb{G}_{A,n}$ with any $|A| < d$ would have no power in the neighborhood of independence.

3.3. *Archimedean copulas.* Following Genest and MacKay [11] and Nelson [20], a copula is called *Archimedean* whenever it can be expressed in the form $C(u) = \phi^{-1}\{\phi(u_1) + \cdots + \phi(u_d)\}$, for some generator $\phi : (0, 1] \to [0, \infty)$. To insure that $C$ is a copula, it suffices that $\phi(1) = 0$ and $(-1)^j d^j \phi^{-1}(t)/dt^j > 0$ for every $j \in \mathcal{S}_d$.

PROPOSITION 3.2. *Let $\mathcal{C}$ be a parametric family of Archimedean copulas with generator $\phi_\theta$ such that $\phi_\theta(t) \to -\log t$ and $\phi'_\theta(t) \to -1/t$ as $\theta \to \theta_0$. Further, assume that $\mathcal{C}$ satisfies assumptions* (i) *and* (ii) *stated in Section* 2. *Then*

$$\frac{\dot{C}_{\theta_0}(u)}{C_{\theta_0}(u)} = \dot{\phi}_{\theta_0}\{C_{\theta_0}(u)\} - \sum_{j=1}^d \dot{\phi}_{\theta_0}(u_j) \tag{3.1}$$

*and*

$$\mu_A(u) = C_{\theta_0}(u^A) \sum_{B \subset A} (-1)^{|A \setminus B|} \dot{\phi}_{\theta_0}\{C_{\theta_0}(u^B)\}. \tag{3.2}$$

PROOF. Noting that $\phi_\theta\{C_\theta(u)\} = \phi_\theta(u_1) + \cdots + \phi_\theta(u_d)$ and applying the chain rule, one finds $\dot{\phi}_\theta\{C_\theta(u)\} + \dot{C}_\theta(u)\phi'_\theta\{C_\theta(u)\} = \dot{\phi}_\theta(u_1) + \cdots + \dot{\phi}_\theta(u_d)$. Equation (3.1) follows by taking the limit as $\theta \to \theta_0$. As for (3.2), it follows from substitution of $\dot{C}_{\theta_0}$ into the formula for $\mathcal{M}_A$, combined with the fact that when $|A| \geq 2$, the Möbius inversion formula yields

$$\sum_{B \subset A} \sum_{j \in B} (-1)^{|A \setminus B|} \dot{\phi}_{\theta_0}(u_j) = \mathbf{1}(|A| = 1) \prod_{j \in A} \dot{\phi}_{\theta_0}(u_j) = 0. \qquad \square$$

EXAMPLE 3.1. Assumptions (i) and (ii) can easily be checked for Frank's family of $d$-variate copulas generated by $\phi_\theta(t) = \log\{(e^{-\theta} - 1)/(e^{-\theta t} - 1)\}$ for $t \in (0, 1]$ with $\theta \in [\ell_d, \infty)$, where $-\infty = \ell_2 < \ell_3 < \cdots < \ell_\infty = 0$. Here, $\theta_0 = 0$ corresponds to independence and $\dot{\phi}_{\theta_0}(t) = (t-1)/2$. Proposition 3.2 and the Möbius inversion formula together yield

$$\dot{C}_{\theta_0}(u) = \frac{1}{2} C_{\theta_0}(u) \left\{ d - 1 + C_{\theta_0}(u) - \sum_{j=1}^d u_j \right\}$$

and

$$\mu_A(u) = \frac{1}{2} C_{\theta_0}(u^A) \sum_{B \subset A} (-1)^{|A \setminus B|} \{C_{\theta_0}(u^B) - 1\}$$

$$= \frac{1}{2} C_{\theta_0}(u^A) \sum_{B \subset A} (-1)^{|A \setminus B|} C_{\theta_0}(u^B) = \frac{1}{2} \prod_{j \in A} u_j(u_j - 1).$$

Here, $\mu_A \neq 0$ for every $A \subset \mathcal{S}_d$ with $|A| \geq 2$.

REMARK 3.1. The $d$-variate Ali–Mikhail–Haq Archimedean system of copulas with $\theta \in [0, 1]$ is generated by $\phi_\theta(t) = (1-\theta)^{-1} \log\{\theta + (1-\theta)/t\}$ for $t \in (0, 1]$ and $\theta_0 = 0$ corresponds to independence. For this system, one gets $\dot\phi_{\theta_0}(t) = t - 1 - \log t$ and hence $\mu_A(u) = \prod_{j \in A} u_j(u_j - 1)$, which is the same as for Frank's family, up to a multiplicative constant. A similar conclusion holds for another version of the Farlie–Gumbel–Morgenstern copula, namely

$$C_\theta(u) = C_{\theta_0}(u) + \theta C_{\theta_0}(u) \left\{ d - 1 + C_{\theta_0}(u) - \sum_{j=1}^d u_j \right\}.$$

EXAMPLE 3.2. Assumptions (i) and (ii) can also be easily verified for Clayton's $d$-variate family of copulas whose generator, defined for all $\theta \in [0, \infty)$, is given by $\phi_\theta(t) = (t^{-\theta} - 1)/\theta$ for $t \in (0, 1]$. Note that independence corresponds to the value $\theta_0 = 0$. One easily finds that $\dot\phi_{\theta_0}(t) = (\log t)^2/2$. Calling on Proposition 3.2 and the Möbius inversion formula, one obtains

$$\dot C_{\theta_0}(u) = C_{\theta_0}(u) \sum_{j<k} \log u_j \log u_k$$

and

$$\mu_A(u) = \mathbf{1}(|A| = 2) \prod_{j \in A} u_j \log u_j.$$

REMARK 3.2. The dependence structure induced by Clayton's copula is also known as the *gamma frailty model*. It may come somewhat as a surprise that in this case, $\mu_A(u) = 0$ unless $|A| = 2$. In other words, Clayton's copula shares with the multivariate normal model the property that tests of independence based on terms $\mathbb{G}_{A,n}$ of the Möbius decomposition with $|A| > 2$ would have no power whatsoever in the neighborhood of independence.

The Gumbel–Barnett system of copulas provides another example of this curious phenomenon. Copulas in this class are generated by $\phi_\theta(t) = \log(1 - \theta \log t)/\theta$ for $t \in (0, 1]$ with $\theta \in [0, 1]$ and $\theta_0 = 0$ corresponding to independence. A simple calculation shows that $\dot\phi_{\theta_0}(t) = -(\log t)^2/2$, and hence the formulas are the same as for Clayton's copula, up to a change in sign.

EXAMPLE 3.3. Assumptions (i) and (ii) can also easily be verified for the so-called Gumbel–Hougaard family of copulas whose generator is defined for $\theta \in [0, 1)$ by $\phi_\theta(t) = |\log t|^{1/(1-\theta)}$ for $t \in (0, 1]$. In that case, one finds $\dot\phi_{\theta_0}(t) = -(\log t) \log(\log 1/t)$. Thus, Proposition 3.2 entails that

$$\dot C_{\theta_0}(u) = C_{\theta_0}(u) \left\{ -\sum_{j=1}^d \log u_j \log\left( \sum_{k=1}^d \log u_k / \log u_j \right) \right\}$$

and

$$\mu_A(u) = -C_{\theta_0}(u^A) \sum_{B \subset A} (-1)^{|A \setminus B|} \left( \sum_{j \in B} \log u_j \right) \log\left( -\sum_{j \in B} \log u_j \right)$$

for all $u \in (0,1)^d$. Here, again, $\mu_A \neq 0$ for every $A \subset \mathcal{S}_d$ with $|A| \geq 2$.

**4. Limiting distributions of Cramér–von Mises functionals.** In the absence of information about the marginal distributions of a multivariate population, a valid testing procedure for independence should be based on some version of the empirical copula process $\mathbb{C}_n$. To improve convergence and reduce bias in finite samples, a centered version of $\mathbb{C}_n$ will be used in the sequel. The latter is defined by

$$\tilde{\mathbb{C}}_n(u) = \sqrt{n}\left\{ C_n(u) - \prod_{j=1}^d U_n(u_j) \right\},$$

where $U_n$ is the cumulative distribution function of a uniformly distributed random variable on the set $\{1/n, \ldots, n/n\}$. It is clear that $\tilde{\mathbb{C}}_n$ and $\mathbb{C}_n$ have the same limiting behavior, so the asymptotic results of Section 2 also apply to $\tilde{\mathbb{C}}_n$ and henceforth to $\tilde{\mathbb{G}}_{A,n} = \mathcal{M}_A(\tilde{\mathbb{C}}_n)$ for $A \subset \mathcal{S}_d$.

A natural way to test for independence is to consider a global measure of discrepancy computed from $\tilde{\mathbb{C}}_n$ or from a combination of distances computed for each of the $\tilde{\mathbb{G}}_{A,n}$ taken individually. Obvious candidates are the Kolmogorov–Smirnov statistics

$$S_n = \sup_{u \in [0,1]^d} |\tilde{\mathbb{C}}_n(u)| \quad \text{and} \quad S_{A,n} = \sup_{u \in [0,1]^d} |\tilde{\mathbb{G}}_{A,n}(u)|$$

and the Cramér–von Mises functionals

$$B_n = \int_{(0,1)^d} \{\tilde{\mathbb{C}}_n(u)\}^2 \, du \quad \text{and} \quad B_{A,n} = \int_{(0,1)^d} \{\tilde{\mathbb{G}}_{A,n}(u)\}^2 \, du.$$

In the sequel, however, attention will be limited to $B_n$ and $B_{A,n}$. Note that these statistics can be expressed as functions of the ranks through

$$B_n = \frac{1}{n} \sum_{i=1}^n \sum_{j=1}^n \prod_{k=1}^d \left( 1 - \frac{R_{ik} \vee R_{jk}}{n} \right)$$

$$- 2 \sum_{i=1}^n \prod_{k=1}^d \left\{ \frac{n(n-1) - R_{ik}(R_{ik} - 1)}{2n^2} \right\}$$

$$+ n \left\{ \frac{(n-1)(2n-1)}{6n^2} \right\}^d$$

and

$$B_{A,n} = \frac{1}{n} \sum_{i=1}^{n} \sum_{j=1}^{n} \prod_{k \in A} D_n(R_{ik}, R_{jk}),$$

where

$$D_n(s,t) = \frac{(n+1)(2n+1)}{6n^2} + \frac{s(s-1)}{2n^2} + \frac{t(t-1)}{2n^2} - \frac{s \vee t}{n}.$$

4.1. *Asymptotics for $B_n$.* Since no explicit Karhunen–Loève expansion for $\mathbb{C}$ is available when $d > 2$, even under $P_n$, the asymptotic null distribution of $B_n$ cannot be computed analytically. This is due to the unwieldy form of the covariance function $\Lambda$.

If $\mathcal{C}$ is a family of copulas whose members satisfy assumptions (i) and (ii), one has from Proposition 2.1 that $B_n$ converges in law under $Q_n$ to

(4.1) $$\mathbb{B} = \int_{(0,1)^d} \{\mathbb{C}(u) + \delta \dot{C}_{\theta_0}(u)\}^2 \, du,$$

which is a Cramér–von Mises functional of a Gaussian process with mean $\delta \dot{C}_{\theta_0}$ and covariance function $\Lambda$. In order to approximate the distribution of $\mathbb{B}$, a procedure due to Deheuvels and Martynov [7] is adopted here. Specifically, the proposed approximation is

(4.2) $$\tilde{\mathbb{B}} = \frac{1}{m} \sum_{i=1}^{m} \xi_i^2,$$

where $\xi = \Delta(U) + V(U)Z$ is an $m$-variate vector constructed from two independent vectors $U$ and $Z$ whose components are independent and $\mathcal{U}(0,1)$ and $\mathcal{N}(0,1)$, respectively. Here, $\Delta(u_1, \ldots, u_m) = \delta(\dot{C}_{\theta_0}(u_1), \ldots, \dot{C}_{\theta_0}(u_m))^\top$ and $V(u_1, \ldots, u_m)$ is the Cholesky decomposition of the covariance matrix $\Sigma$ with components $\Sigma_{jk} = \Lambda(u_j, u_k)$ for $j, k \in \mathcal{S}_m$. Deheuvels and Martynov show that the variance of the approximation error is $O(1/m)$.

4.2. *Asymptotics for $B_{A,n}$.* In view of the fact that the covariance structure of $\mathbb{G}_A$ is a product of covariance functions associated with Brownian bridges, it follows from standard theory [22, p. 213] that under $P_n$, the limiting process $\mathbb{G}_A$ admits the representation

(4.3) $$\mathbb{G}_A(u) = \sum_{\gamma \in \mathbb{N}^{|A|}} \sqrt{\lambda_\gamma} Z_\gamma f_\gamma(u), \qquad \gamma = (\gamma_j)_{j \in A},$$

where the $Z_\gamma$ are independent $\mathcal{N}(0,1)$ random variables and

$$\lambda_\gamma = \prod_{j \in A} (\pi \gamma_j)^{-2}, \qquad f_\gamma(u) = \prod_{j \in A} \sqrt{2} \sin(\gamma_j \pi u_j), \qquad u \in [0,1]^d.$$

An idea of Deheuvels [6], later exploited by Genest and Rémillard [13], is to base a test of independence on some combination of the asymptotically independent statistics $B_{A,n}$. In view of Corollary 2.2, the limiting distribution of $B_{A,n}$ under $Q_n$ is given by

$$\mathbb{B}_A = \int_{(0,1)^d} \{\mathbb{G}_A(u) + \delta\mu_A(u)\}^2 \, du.$$

As already shown by Deheuvels [6], formula (4.3) implies that

$$\mathbb{Q}_A = \int_{(0,1)^d} \{\mathbb{G}_A(u)\}^2 \, du = \sum_{\gamma \in \mathbb{N}^{|A|}} \lambda_\gamma Z_\gamma^2, \qquad \lambda_\gamma = \prod_{j \in A} \lambda_{\gamma_j}.$$

It follows that

$$\mathbb{B}_A = \mathbb{Q}_A + 2\delta \sum_{\gamma \in \mathbb{N}^{|A|}} \lambda_\gamma I_{\gamma,A} Z_\gamma + \delta^2 I_A,$$

where

$$I_A = \int_{(0,1)^d} \{\mu_A(u)\}^2 \, du \quad \text{and} \quad I_{\gamma,A} = \frac{1}{\sqrt{\lambda_\gamma}} \int_{(0,1)^d} \mu_A(u) f_\gamma(u) \, du.$$

The limiting distribution of $B_{A,n}$ under $Q_n$ is given below. The result follows from direct substitution of Parseval's identity $I_A = \sum_{\gamma \in \mathbb{N}^{|A|}} \lambda_\gamma I_{\gamma,A}^2$ into the integral representation of $\mathbb{B}_A$.

PROPOSITION 4.1. *If $\mathcal{C}$ is a family of copulas whose members satisfy assumptions* (i) *and* (ii), *then the asymptotic distribution of $B_{A,n}$ under $Q_n$ is given by the weighted sum*

$$\mathbb{B}_A = \sum_{\gamma \in \mathbb{N}^{|A|}} \lambda_\gamma (Z_\gamma + \delta I_{\gamma,A})^2.$$

Note that there also exists a representation like (4.3) for $\mathbb{C}$, but its weights are unknown. Nevertheless, their sum is finite. Once again calling on Parseval's identity and using (4.1), one obtains

$$\mathbb{B} = \sum_{\gamma \in \mathbb{N}^d} \tilde{\lambda}_\gamma (\tilde{Z}_\gamma + \delta \tilde{I}_\gamma)^2,$$

where

$$\sum_{\gamma \in \mathbb{N}^d} \tilde{\lambda}_\gamma = \int_{(0,1)^d} \int_{(0,1)^d} \Gamma(u, u') \, du' \, du < \infty$$

and

$$\sum_{\gamma \in \mathbb{N}^d} \tilde{\lambda}_\gamma \tilde{I}_\gamma^2 = \int_{(0,1)^d} \{\dot{C}_{\theta_0}(u)\}^2 \, du < \infty.$$

The unknown quantities $\tilde{\lambda}_\gamma$ and $\tilde{I}_\gamma$ could be approximated numerically.

For many systems of distributions, the drift term has the simple form $\mu_A(u) = \prod_{j \in A} \mu(u_j)$. In that case, it is easy to see that

$$(4.4) \quad I_{\gamma,A} = \frac{2^{|A|/2}}{\sqrt{\lambda_\gamma}} \prod_{j \in A} f(\gamma_j), \qquad \text{where } f(k) = \int_0^1 \mu(u) \sin(k\pi u)\, du.$$

EXAMPLE 4.1 (Equicorrelated Gaussian copulas). For this model, (4.4) and results from Section 3.1 together imply that $I_{\gamma,A} = \mathbf{1}(|A| = 2) 2\pi^2 \prod_{j \in A} \gamma_j g(\gamma_j)$, where $g(k) = \int \varphi\{\Phi^{-1}(u)\} \sin(k\pi u)\, du$.

EXAMPLE 4.2 (Farlie–Gumbel–Morgenstern copulas). For this model, results from Section 3.2 imply that $I_{\gamma,A} = (-1)^d \mathbf{1}(A = \mathcal{S}_d) 2^{5d/2} \lambda_\gamma$ if $\gamma_1$, ..., $\gamma_d$ are all odd and $I_{\gamma,A} = 0$ otherwise.

EXAMPLE 4.3 (Frank and Ali–Mikhail–Haq copulas). For these models, calculations based on material from Example 3.1 lead to $I_{\gamma,A} = (-1)^{|A|} 2^{5|A|/2 - 1} \lambda_\gamma$ if $\gamma_j$ is odd for every $j \in A$ and $I_{\gamma,A} = 0$ otherwise.

EXAMPLE 4.4 (Clayton and Gumbel–Barnett distributions). For these models, the observations already made in Example 3.3 and Remark 3.2 yield $I_{\gamma,A} = \mathbf{1}(|A| = 2) 2\pi^{-2} \prod_{j \in A} SI(\gamma_j \pi)/\gamma_j$, where $SI(x) = \int_0^x t^{-1} \sin(t)\, dt$.

4.3. *Combination of independent statistics.* A simple test of independence consists of rejecting the null hypothesis whenever the observed value of $B_n$ exceeds the $(1 - \alpha)$-percentile of its asymptotic distribution under $P_n$. The critical value $q_B(\alpha)$ such that $\mathbf{P}\{B_n > q_B(\alpha) | P_n\} \to \alpha$ as $n \to \infty$ is easily approximated using formula (4.2) with $\delta = 0$. For $\alpha = 0.05$, these values are to be found in Table 1.

However, potentially more efficient methods could be based on combinations of the asymptotically independent $B_{A,n}$'s with $A \in \mathcal{S}_d$. Four such procedures are considered here. The first is inspired by Ghoudi, Kulperger and Rémillard [14], while the second and third were studied by Genest and Rémillard [13].

4.3.1. *Linear combination rule.* Base the test on

$$L_n = \sum_{|A| > 1} B_{A,n}.$$

Table 1 gives approximate values for $q_L(\alpha)$ such that $\mathbf{P}\{L_n > q_L(\alpha) | P_n\} \to \alpha$ as $n \to \infty$.

TABLE 1
*Approximate critical values of the tests based on $B_n$, $L_n$, $W_n$, $M_n$ and $T_n$, at the 5% level*

|  | $d=3$ | $d=4$ | $d=5$ |
|---|---|---|---|
| $q_B(0.05)$ | 0.05669 | 0.04124 | 0.02549 |
| $q_L(0.05)$ | 0.14045 | 0.26002 | 0.42046 |
| $q_W(0.05)$ | 18.50403 | 50.54507 | 134.38756 |
| $q_2(\alpha')$ | 0.08518 | 0.10126 | 0.12731 |
| $q_3(\alpha')$ | 0.01006 | 0.01186 | 0.01326 |
| $q_4(\alpha')$ | – | 0.00159 | 0.00175 |
| $q_5(\alpha')$ | – | – | 0.00022 |
| $q_T(0.05)$ | 15.31231 | 35.09049 | 71.00888 |

4.3.2. *Dependogram method.* Base the test on

$$M_n = \max_{|A|>1}\{B_{A,n}/q_{|A|}(\alpha')\},$$

where $q_{|A|}(\alpha)$ is such that $\mathbf{P}\{B_{A,n} > q_{|A|}(\alpha)|P_n\} \to \alpha$ as $n \to \infty$ and where

(4.5) $$\alpha' = 1 - (1-\alpha)^{1/(2^d-d-1)}$$

is chosen so that by performing each test at level $\alpha'$, the global level of the procedure is $\alpha$. The first few simulated critical values $q_{|A|}$ are reported in Table 1. Note that $\mathbf{P}(M_n > 1|P_n) \to \alpha$ as $n \to \infty$.

4.3.3. *Fisher's approach.* Base the test on

$$T_n = -2 \sum_{|A|>1} \log\{1 - F_{A,n}(B_{A,n})\},$$

where $F_{A,n}(x) = \mathbf{P}(B_{A,n} \le x)$. Under $P_n$, $T_n$ converges in law to

$$-2 \sum_{|A|>1} \log\{1 - q_{|A|}^{-1}(\mathbb{Q}_A)\},$$

which is chi-square with $2(2^d - d - 1)$ degrees of freedom. Hence, the critical value of a test based on $T_n$ is given by

$$q_T(\alpha) = K^{-1}(1-\alpha), \qquad \text{where } K(x) = \mathbf{P}\{\chi^2_{2(2^d-d-1)} \le x\}.$$

Table 1 reports the asymptotic values of $q_T(0.05)$ for $d \in \{3,4,5\}$, as per Genest and Rémillard [13].

4.3.4. *Weighted linear combination rule.* Base the test on

$$W_n = \sum_{|A|>1} \pi^{2|A|} B_{A,n}.$$

Table 1 gives approximate values for $q_W(\alpha)$ such that $\mathbf{P}\{W_n > q_W(\alpha)|P_n\} \to \alpha$ as $n \to \infty$. This procedure is inspired by Proposition A.1 since in view of the latter, $-2\log\{1 - F_{A,n}(B_{A,n})\} \approx \pi^{2|A|} B_{A,n}$ whenever $B_{A,n}$ is large. Therefore, $W_n$ should approximate $T_n$ well under fixed alternatives $C_\theta$ because in that case, $B_{A,n}/n$ tends to a positive constant.

**5. Comparison of local power functions.** One way to compare competing test procedures for independence is to evaluate their asymptotic power function in a neighborhood of $\theta = \theta_0$, that is, under copula alternatives $C_{\theta_n}$ that form a contiguous sequence $(Q_n)$. Specifically, let $S_n$ be some statistic for independence with asymptotic critical value $q_S(\alpha)$. The associated local power function is defined as

$$\beta_S(\alpha, \delta) = \lim_{n \to \infty} \mathbf{P}\{S_n > q_S(\alpha)|Q_n\}.$$

Analytic expressions for this function are given in Section 5.1 for the combined test statistics $L_n$, $M_n$ and $W_n$. No similar form could be obtained for $B_n$ or for Fisher's procedure $T_n$; however, see Remark 5.1. For numerical comparisons, see Sections 5.2 and 5.3.

5.1. *Analytic expressions for local power functions.* A result that will prove useful in the sequel is the formula of Gil-Pelaez [16], which says that if $X$ is a random variable with characteristic function $\phi_X$, then

$$\mathbf{P}(X > x) = \frac{1}{2} + \frac{1}{2\pi} \int_{-\infty}^{\infty} \mathrm{Im}\{t^{-1} e^{-ixt} \phi_X(t)\} dt,$$

where $\mathrm{Im}(z)$ denotes the imaginary part of any complex number $z$.

To use this identity in the present context, let

$$\hat{\eta}(t, \delta) = (1 - 2it)^{-1/2} \exp\left(\frac{i\delta^2 t}{1 - 2it}\right)$$
$$= (1 + 4t^2)^{-1/4} e^{-2t^2\delta^2/(1+4t^2)} e^{i\arctan(2t)/2 + it\delta^2/(1+4t^2)}$$

be the characteristic function of a noncentral chi-square variable $(Z + \delta)^2$, where $Z \sim \mathcal{N}(0, 1)$. From Proposition 4.1, it follows that

$$\phi_{\mathbb{B}_A}(t, \delta) = \mathrm{E}(e^{it\mathbb{B}_A}) = \prod_{\gamma \in \mathbb{N}^{|A|}} \hat{\eta}(\lambda_\gamma t, \delta I_{\gamma, A})$$
$$= \xi_A(t) e^{-2\delta^2 t^2 \kappa_{A,1}(t)} e^{i\kappa_{A,2}(t) + i\delta^2 \kappa_{A,3}(t)},$$

where

$$\xi_A(t) = \prod_{\gamma \in \mathbb{N}^{|A|}} (1 + 4t^2 \lambda_\gamma^2)^{-1/4}, \qquad \kappa_{A,1}(t) = \sum_{\gamma \in \mathbb{N}^{|A|}} \lambda_\gamma^2 I_{\gamma,A}^2 / (1 + 4t^2 \lambda_\gamma^2),$$

$$\kappa_{A,2}(t) = \frac{1}{2} \sum_{\gamma \in \mathbb{N}^{|A|}} \arctan(2t\lambda_\gamma), \qquad \kappa_{A,3}(t) = t \sum_{\gamma \in \mathbb{N}^{|A|}} \lambda_\gamma I_{\gamma,A}^2 / (1 + 4t^2 \lambda_\gamma^2).$$

According to Proposition 4.1, the limiting distribution of $L_n$ under $Q_n$ is $L = \sum_{|A|>1} \mathbb{B}_A$, so its asymptotic characteristic function is given by

$$\phi_L(t,\delta) = \prod_{|A|>1} \phi_{\mathbb{B}_A}(t,\delta) = \xi(t) e^{-2\delta^2 t^2 \kappa_1(t)} e^{i\kappa_2(t) + i\delta^2 \kappa_3(t)},$$

where $\xi(t) = \prod_{|A|>1} \xi_A(t)$ and $\kappa_i(t) = \sum_{|A|>1} \kappa_{A,i}(t)$ for $i \in \{1,2,3\}$.

An application of the Gil–Pelaez formula then yields the following result:

PROPOSITION 5.1. *If $\mathcal{C}$ is a family of copulas whose members satisfy assumptions* (i) *and* (ii), *then under $Q_n$, $\beta_L(\alpha, \delta) = \mathbf{P}\{L > q_L(\alpha)\}$, where*

$$\mathbf{P}(L > x) = \frac{1}{2} + \frac{1}{\pi} \int_0^\infty \frac{\sin\{\kappa(x,t)\}}{t\zeta(t)} dt$$

with

$$\kappa(x,t) = -\frac{xt}{2} + \frac{1}{2} \sum_{|A|>1} \sum_{\gamma \in \mathbb{N}^{|A|}} \left\{ \arctan(\lambda_\gamma t) + \delta^2 \frac{\lambda_\gamma I_{\gamma,A}^2 t}{1 + \lambda_\gamma^2 t^2} \right\}$$

and

$$\zeta(t) = \exp\left( \frac{\delta^2 t^2}{2} \sum_{|A|>1} \sum_{\gamma \in \mathbb{N}^{|A|}} \frac{\lambda_\gamma^2 I_{\gamma,A}^2}{1 + \lambda_\gamma^2 t^2} \right) \prod_{|A|>1} \prod_{\gamma \in \mathbb{N}^{|A|}} (1 + \lambda_\gamma^2 t^2)^{1/4}.$$

A similar result for the test based on $M_n$ follows from the fact that

$$\lim_{n \to \infty} \mathbf{P}(M_n > 1) = 1 - \prod_{|A|>1} \mathbf{P}\{\mathbb{B}_A \leq q_{|A|}(\alpha')\}$$

(5.1)
$$= 1 - \prod_{|A|>1} \{1 - \beta_A(\alpha', \delta)\}$$

with $\alpha'$ defined as in (4.5).

PROPOSITION 5.2. *Let $M$ be the limit in distribution of $M_n$ under $Q_n$ and define $\beta_M$ to be its associated local power function. If $\mathcal{C}$ is a family*

ALE OF TESTS OF INDEPENDENCE                                    17of copulas whose members satisfy assumptions (i) and (ii), then under $Q_n$, $\beta_M(\alpha, \delta)$ is given by (5.1) and for all $x > 0$,

$$\mathbf{P}(\mathbb{B}_A > x) = \frac{1}{2} + \frac{1}{\pi} \int_0^\infty \frac{\sin\{\kappa_A(x,t)\}}{t\zeta_A(t)} dt$$

with

$$\kappa_A(x,t) = -\frac{xt}{2} + \frac{1}{2} \sum_{\gamma \in \mathbb{N}^{|A|}} \left\{ \arctan(\lambda_\gamma t) + \delta^2 \frac{\lambda_\gamma I_{\gamma,A}^2 t}{1 + \lambda_\gamma^2 t^2} \right\}$$

and

$$\zeta_A(t) = \exp\left(\frac{\delta^2 t^2}{2} \sum_{\gamma \in \mathbb{N}^{|A|}} \frac{\lambda_\gamma^2 I_{\gamma,A}^2}{1 + \lambda_\gamma^2 t^2}\right) \prod_{\gamma \in \mathbb{N}^{|A|}} (1 + \lambda_\gamma^2 t^2)^{1/4}.$$

Next, one can see that the limiting distribution of $W_n$ under $Q_n$ is $W = \sum_{|A|>1} \pi^{2|A|} \mathbb{B}_A$, so its asymptotic characteristic function is given by

$$\prod_{|A|>1} \phi_{\mathbb{B}_A}(t\pi^{2|A|}, \delta) = \xi_W(t) e^{-2\delta^2 t^2 \pi^{4|A|} \kappa_{1,W}(t)} e^{i\kappa_{2,W}(t) + i\delta^2 \kappa_{3,W}(t)},$$

where for $i \in \{1, 2, 3\}$,

$$\xi_W(t) = \prod_{|A|>1} \xi_A(t\pi^{2|A|}) \quad \text{and} \quad \kappa_{i,W}(t) = \sum_{|A|>1} \kappa_{A,i}(t\pi^{2|A|}).$$

The following proposition gives the local power function of $W$:

PROPOSITION 5.3. *If $\mathcal{C}$ is a family of copulas whose members satisfy assumptions* (i) *and* (ii)*, then under $Q_n$, $\beta_W(\alpha, \delta) = \mathbf{P}\{W > q_W(\alpha)\}$, where*

$$\mathbf{P}(W > x) = \frac{1}{2} + \frac{1}{\pi} \int_0^\infty \frac{\sin\{\kappa_W(x,t)\}}{t\zeta_W(t)} dt$$

with

$$\kappa_W(x,t) = -\frac{xt}{2} + \frac{1}{2} \sum_{|A|>1} \sum_{\gamma \in \mathbb{N}^{|A|}} \left\{ \arctan(\lambda_\gamma \pi^{2|A|} t) + \delta^2 \frac{\lambda_\gamma \pi^{2|A|} I_{\gamma,A}^2 t}{1 + \lambda_\gamma^2 \pi^{4|A|} t^2} \right\}$$

and

$$\zeta(t) = \exp\left(\frac{\delta^2 t^2}{2} \sum_{|A|>1} \sum_{\gamma \in \mathbb{N}^{|A|}} \frac{\lambda_\gamma^2 \pi^{4|A|} I_{\gamma,A}^2}{1 + \lambda_\gamma^2 \pi^{4|A|} t^2}\right) \prod_{|A|>1} \prod_{\gamma \in \mathbb{N}^{|A|}} (1 + \lambda_\gamma^2 \pi^{4|A|} t^2)^{1/4}.$$

REMARK 5.1. For Fisher's test, no explicit representation for $\beta_T(\alpha, \delta)$ seems possible. Under $Q_n$, however, $T_n$ converges in distribution to

$$-2 \sum_{|A|>1} \log\{1 - q_{|A|}^{-1}(\mathbb{B}_A)\},$$

where

$$q_{|A|}^{-1}(x) = \frac{1}{2} - \frac{1}{\pi} \int_0^\infty \frac{\sin\{\kappa_0(x,t)\}}{t\zeta_0(t)} \, dt$$

with

$$\kappa_0(x,t) = -\frac{xt}{2} + \frac{1}{2} \sum_{\gamma \in \mathbb{N}^{|A|}} \arctan(\lambda_\gamma t) \quad \text{and} \quad \zeta_0(t) = \prod_{\gamma \in \mathbb{N}^{|A|}} (1 + \lambda_\gamma^2 t^2)^{1/4}.$$

An approximation is then given by

$$\hat{\beta}_T(\alpha, \delta) = \frac{1}{N} \sum_{i=1}^N \mathbf{1}\{T_i > q_T(\alpha)\},$$

where

$$T_i = -2 \sum_{|A|>1} \log\{1 - q_{|A|}^{-1}(\mathbb{B}_{A,i})\}$$

and for $N$ sufficiently large, $\mathbb{B}_{A,1}, \ldots, \mathbb{B}_{A,N}$ are mutually independent observations from a finite-sum approximation of $\mathbb{B}_A$.

5.2. *Power comparisons for statistics involving all $A \subset \mathcal{S}_d$ with $|A| > 1$.*
Figure 1 graphically compares the power of five different tests for trivariate independence based on statistics $B_n$, $L_n$, $M_n$, $T_n$ and $W_n$. These comparisons were carried out at the 5% level for four parametric classes of three-dimensional copula alternatives considered in Section 3, namely (a) the equicorrelated Gaussian copula, (b) the Farlie–Gumbel–Morgensten copula, (c) the Frank or Ali–Mikhail–Haq copula and (d) the Clayton or Gumbel–Barnett copula. In view of the considerations made in Sections 3 and 4, the choice of representative within each class is irrelevant for the following asymptotic local power comparisons.

The local power curves corresponding to $L_n$, $M_n$ and $W_n$ were evaluated by means of numerical integration of the formulas given in Propositions 5.1, 5.2 and 5.3. In each case, the integral was computed by the trapezoidal rule on a domain of the form $[0, K]$ for suitably large $K$; all infinite sums were truncated to 40 terms in each index. This was sufficient to ensure numerical accuracy. As for the local power curves associated with $B_n$ and $T_n$, they were obtained via Monte Carlo simulation, using 10,000 repetitions for each point on the curve.

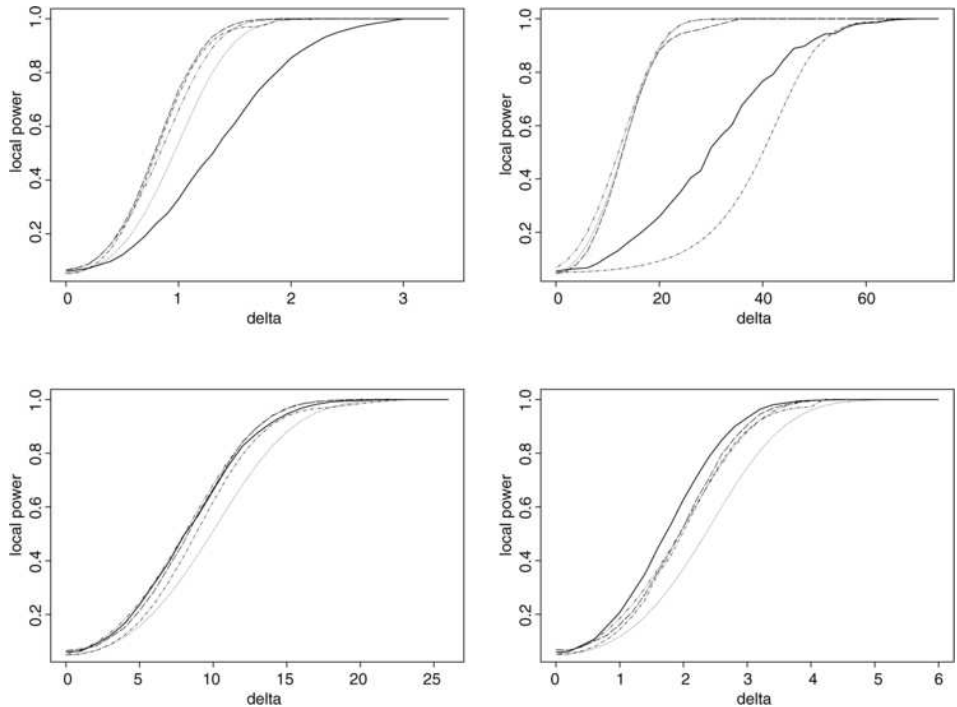

Fig. 1. *Comparative local power of $B_n$ (—), $T_n$ (- -), $M_n$ (· · ·), $L_n$ (-·-) and $W_n$ (-··-) for four different classes of trivariate copula alternatives: the equicorrelated Gaussian (upper left panel), the Farlie–Gumbel–Morgenstern (upper right), Frank/Ali–Mikhail–Haq (lower left) and the Clayton/Gumbel–Barnett copula (lower right).*

Looking at Figure 1, one can see that the test based on $B_n$ is best for the Clayton/Gumbel–Barnett models in the neighborhood of independence. It is also close to optimal for the Frank/Ali–Mikhail–Haq models, but its behavior is much less satisfactory for the Farlie–Gumbel–Morgenstern copula and especially for the equicorrelated Gaussian structure.

Putting $B_n$ aside, the smallest local power is yielded by $M_n$ for all structures considered except the Farlie–Gumbel–Morgenstern. In the latter model, $L_n$ is the least powerful locally. This results from the fact that this procedure equally weights the four $B_{A,n}$ statistics of size $|A| = 2$ or $3$, whereas only the latter has any power in detecting dependence, as observed in Section 3.2. In contrast, note the excellent performance of $W_n$, which weights the $B_{A,n}$ proportionally to $\pi^{2|A|}$.

Given that $W_n$ is an approximation of $T_n$, as seen in Section 4.3.4, it is not surprising that they behave similarly in all cases. While Fisher's procedure $T_n$ ranks first in the trivariate normal model, there does not seem to be much difference between $T_n$ and $W_n$ in the other models considered.

TABLE 2
*Approximate critical values of the tests based on $L_{n,2}$, $M_{n,2}$ and $T_{n,2}$, at the 5% level*

|  | $d=3$ | $d=4$ | $d=5$ |
|---|---|---|---|
| $q_{L_2}(0.05)$ | 0.13562 | 0.23917 | 0.36848 |
| $q_2(\alpha'')$ | 0.07479 | 0.09020 | 0.09714 |
| $q_{T_2}(0.05)$ | 12.20343 | 20.32429 | 29.63573 |

5.3. *Power comparisons for statistics based on sets $A$ with $|A|=2$.* In Section 3, it was seen that for Clayton/Gumbel–Barnett and equicorrelated Gaussian contiguous alternatives, $\mu_A = 0$ when $|A| > 2$. As a consequence, the asymptotic distribution of $B_{A,n}$ under $Q_n$ is then the same as that under $P_n$. A loss in efficiency may thus be expected to occur when testing for independence within these models whenever a test statistic combines all possible $B_{A,n}$, rather than only those for which $|A|=2$.

For such alternatives, possibly more efficient tests could be based on

$$L_{n,2} = \sum_{|A|=2} B_{A,n}, \qquad M_{n,2} = \max_{|A|=2} \left\{ \frac{B_{A,n}}{q_2(\alpha'')} \right\}$$

and

$$T_{n,2} = \sum_{|A|=2} \log\{1 - F_{A,n}(B_{A,n})\},$$

where $\alpha'' = 1 - (1-\alpha)^{2/\{d(d-1)\}}$ and $q_2(\alpha'')$ is given in Table 2. Analytic expressions for the asymptotic local power curves of these statistics can be derived from straightforward adaptations of Propositions 5.1, 5.2 and 5.3. Table 2 gives the critical values, at the 5% level, of the tests based on these statistics, based on Monte Carlo simulations. These values were used in Figure 2 to compare the local power curves of these statistics to their analogous versions based on all $|A| > 1$.

In the case $d = 3$, only the statistic $B_{A,n}$ with $|A| = 3$ has no local power in testing for independence in the equicorrelated Gaussian model. The loss of power caused by its inclusion is most apparent for the dependogram statistic $M_n$, as illustrated in the left panel of Figure 2. The right panel of the same figure shows that the loss is much less dramatic (if indeed there is any) for $T_n$ compared to $T_{n,2}$. This was perhaps to be expected, considering that Fisher's test based on $T_n$ is probably close to optimal in this context.

**6. Asymptotic efficiency results.** Given two competing test statistics $S_{n,1}$ and $S_{n,2}$, let

$$\beta_{S_i}(\alpha, \delta) = \lim_{n \to \infty} \mathbf{P}\{S_{n,i} > q_{S_i}(\alpha) | Q_n\}$$

ALE OF TESTS OF INDEPENDENCE 21

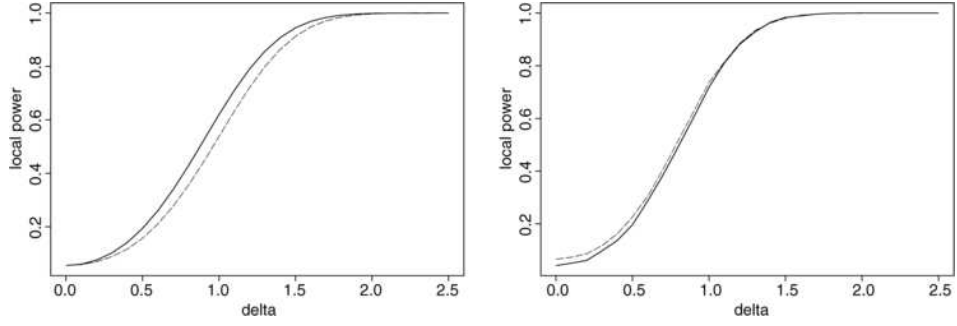

Fig. 2. *Comparison of the local power of two test statistics for independence in the trivariate, equicorrelated Gaussian model. Left panel: $M_{n,2}$ (solid line) versus $M_n$ (broken line). Right panel: $T_{n,2}$ (solid line) versus $T_n$ (broken line).*

and suppose that $\{\beta_{S_i}(\alpha,\delta) - \alpha\}/\delta \to 0$ as $\delta \to 0$ for $i \in \{1,2\}$. Genest, Quessy and Rémillard [12] argue that an appropriate measure of the asymptotic relative efficiency of $S_{n,1}$ with respect to $S_{n,2}$ is then given by

$$(6.1) \quad e_{12}(\alpha) = \lim_{\delta \to 0} \frac{\beta_{S_1}(\alpha,\delta) - \alpha}{\beta_{S_2}(\alpha,\delta) - \alpha} = \frac{\lim_{\delta \to 0}\{\beta_{S_1}(\alpha,\delta) - \alpha\}/\delta^2}{\lim_{\delta \to 0}\{\beta_{S_2}(\alpha,\delta) - \alpha\}/\delta^2},$$

which is a limiting ratio of curvatures of local power functions at $\delta = 0$. Analytic expressions for these curvatures can be obtained for the tests of independence based on $L_n$, $W_n$ and $M_n$. They are given in the next proposition, whose proof depends on Proposition A.2. In the sequel, $\beta_A$ denotes the local power function associated with the statistic $B_{A,n}$.

PROPOSITION 6.1. *Let $h_{w,|A|}$, $h_{w,L}$ and $h_{w,W}$ be the densities of*

$$\mathbb{Q}_A + w\chi_2^2, \quad \sum_{|A|>1} \mathbb{Q}_A + w\chi_2^2 \quad \text{and} \quad \sum_{|A|>1} \pi^{2|A|}\mathbb{Q}_A + w\chi_2^2,$$

*respectively, where in each case, the summands are taken to be independent. Then as $\delta \to 0$, one has*

(i)
$$\{\beta_A(\alpha,\delta) - \beta_A(\alpha,0)\}/\delta^2 \to \sum_{\gamma \in \mathbb{N}^{|A|}} \lambda_\gamma I_{\gamma,A}^2 h_{\lambda_\gamma,|A|}\{q_{|A|}(\alpha)\};$$

(ii)
$$\{\beta_L(\alpha,\delta) - \beta_L(\alpha,0)\}/\delta^2 \to \sum_{|A|>1}\sum_{\gamma \in \mathbb{N}^{|A|}} \lambda_\gamma I_{\gamma,A}^2 h_{\lambda_\gamma,L}\{q_L(\alpha)\};$$

(iii)
$$\{\beta_W(\alpha,\delta) - \beta_W(\alpha,0)\}/\delta^2 \to \sum_{|A|>1}\sum_{\gamma \in \mathbb{N}^{|A|}} \lambda_\gamma \pi^{2|A|} I_{\gamma,A}^2 h_{\lambda_\gamma,W}\{q_W(\alpha)\};$$

(iv) *given $\alpha'$ defined as in* (4.5),
$$\{\beta_M(\alpha,\delta) - \beta_M(\alpha,0)\}/\delta^2 \to \left(\frac{1-\alpha}{1-\alpha'}\right) \sum_{|A|>1} \sum_{\gamma \in \mathbb{N}^{|A|}} \lambda_\gamma I_{\gamma,A}^2 h_{\lambda_\gamma,|A|}\{q_{|A|}(\alpha')\}.$$

PROOF.  Recalling that
$$\mathbb{B}_A = \sum_{\gamma \in \mathbb{N}^{|A|}} \lambda_\gamma (Z_\gamma + \delta I_{\gamma,A})^2, \qquad L = \sum_{|A|>1} \mathbb{B}_A \quad \text{and} \quad W = \sum_{|A|>1} \pi^{2|A|} \mathbb{B}_A,$$

one can derive results (i), (ii) and (iii) easily using Proposition A.2. Now, using (i), one has $\beta_A(\alpha,\delta) = \beta_A(\alpha,0) + \delta^2 \ell_A(\alpha) + o(\delta^2)$, where
$$\ell_A(\alpha) = \sum_{\gamma \in \mathbb{N}^{|A|}} \lambda_\gamma I_{\gamma,A}^2 h_{\lambda_\gamma,|A|}\{q_{|A|}(\alpha)\}.$$

Statement (iv) then follows, on writing
$$\begin{aligned}
\beta_M(\alpha,\delta) - \beta_M(\alpha,0) &= \prod_{|A|>1}\{1 - \beta_A(\alpha',0)\} - \prod_{|A|>1}\{1 - \beta_A(\alpha',\delta)\} \\
&= \delta^2 \sum_{|A|>1} \ell_A(\alpha') \prod_{|D|>1, D \neq A}\{1 - \beta_D(\alpha',0)\} + o(\delta^2) \\
&= \left(\frac{1-\alpha}{1-\alpha'}\right)\delta^2 \sum_{|A|>1} \ell_A(\alpha') + o(\delta^2). \qquad \square
\end{aligned}$$

The above result, combined with equation (6.1), makes it possible to determine the local asymptotic efficiencies of statistics $L_n$, $M_n$ and $W_n$. These values are given in Table 3, along with those corresponding to $L_{n,2}$ and $M_{n,2}$, for which a simple adaptation of Proposition 6.1 must be used.

The results in Table 3 are generally in line with those reported in Section 5, bearing in mind that due to the lack of analytical expressions for their local power curves, the Cramér–von Mises and Fisher test procedures, based respectively on $B_n$ and $T_n$, could not be included in the efficiency comparisons.

The following observations are offered as concluding remarks.

(a) Among trivariate tests of independence based on $L_n$, $M_n$ and $W_n$ which involve all of the $B_{A,n}$, $L_n$ was most efficient and $M_n$ was least efficient in all models considered except the Farlie–Gumbel–Morgenstern class of copulas. In the latter case, $M_n$ was by far the best choice.

(b) Tests based on $L_{n,2}$ or $M_{n,2}$ were totally inefficient when dependence entered through Farlie–Gumbel–Morgenstern trivariate copulas. This was to be expected because in this case, only the term involving $B_{A,n}$ with $|A| = 3$ has a limiting distribution that differs under the null and under the alternative. Tests based on $L_n$ and $W_n$ hardly did any better.

TABLE 3
*Local asymptotic relative efficiency of Cramér–von Mises tests for trivariate independence*

| Alternative trivariate model | Best statistic among $L_n, M_n, W_n, L_{n,2}, M_{n,2}$ for the chosen model | Relative efficiency of the statistic with respect to the best | | | | |
|---|---|---|---|---|---|---|
| | | $L_n$ | $L_{n,2}$ | $M_n$ | $M_{n,2}$ | $W_n$ |
| Equicorrelated Gaussian | $L_{n,2}$ | 98.56 | 100.0 | 44.66 | 88.78 | 79.45 |
| Farlie–Gumbel–Morgenstern | $M_n$ | 3.71 | 0 | 100.0 | 0 | 32.95 |
| Frank Ali–Mikhail–Haq | $L_{n,2}$ | 99.34 | 100.0 | 65.28 | 88.92 | 86.09 |
| Clayton Gumbel–Barnett | $L_{n,2}$ | 98.55 | 100.0 | 43.27 | 87.21 | 79.79 |

(c) For the other models, a loss in efficiency was observed when going from $L_{n,2}$ to $L_n$ and from $M_{n,2}$ to $M_n$. In the equicorrelated Gaussian and Clayton/Gumbel–Barnett cases, this was expected since it was argued earlier that the inclusion of statistics $B_{A,n}$ with $|A| > 2$ is then likely to dilute the power of the overall procedure. An explanation for the same occurrence in the Frank and Ali–Mikhail–Haq models is still lacking.

If nothing else, this study provides a new illustration of the truism that no single procedure could ever be declared best for simultaneously testing for all forms of multivariate dependence. More importantly, however, the results reported herein shed new light on unsuspected commonalities among classes of dependence models that may be superficially perceived as quite different. The accumulation of evidence from this and similar investigations may eventually lead to new typologies for dependence. Needless to say, much remains to be done before this goal can be achieved.

## APPENDIX: AUXILIARY RESULTS

The following result justifies combination procedure 4.3.4.

PROPOSITION A.1. *Suppose that $X = \sum_{k=1}^{\infty} w_k Z_k^2$, where the $Z_k$ are independent $\mathcal{N}(0,1)$ random variables, $w_1 \geq w_2 \geq \cdots$ and $\mathrm{E}(X) = \sum_{k=1}^{\infty} w_k < \infty$. Then $x^{-1} \log \mathbf{P}(X > x) \to -1/(2w_1)$ as $x \to \infty$.*

PROOF. Fix $\alpha < 1/(2w_1)$. Since $\log(1-x) \leq x/(1-x)$, one has

$$\log\{\mathrm{E}(e^{\alpha X})\} = -\frac{1}{2} \sum_{k=1}^{\infty} \log(1 - 2\alpha w_k) \leq \frac{\alpha \mathrm{E}(X)}{1 - 2\alpha w_1} < \infty.$$

Hence, by Markov's inequality one obtains $\mathbf{P}(X > x) \leq e^{-\alpha x}\mathrm{E}(e^{\alpha X})$ and so $\limsup_{x\to\infty} x^{-1}\log\mathbf{P}(X > x) \leq -\alpha$. Letting $\alpha \to 1/(2w_1)$, one obtains $\limsup_{x\to\infty} x^{-1}\log\mathbf{P}(X > x) \leq -1/(2w_1)$. This concludes the proof since in view of large deviation results for Gaussian variables, one also has

$$\liminf_{x\to\infty} x^{-1}\log\mathbf{P}(X > x) \geq \liminf_{x\to\infty} x^{-1}\log\mathbf{P}(w_1 Z_1^2 > x) = -\frac{1}{2w_1}. \qquad \square$$

The following result is instrumental in establishing Proposition 6.1.

PROPOSITION A.2. *Let $(w_k)$ be a positive sequence with $\sum_{k=1}^\infty w_k < \infty$ and $(\mu_k)$ be a sequence such that $\sum_{k=1}^\infty w_k \mu_k < \infty$. Furthermore, let $(Z_k)$ be a sequence of independent $\mathcal{N}(0,1)$ random variables. For any $\delta \geq 0$, set*

$$X_\delta = \sum_{k=1}^\infty w_k (Z_k + \delta\mu_k)^2.$$

*Then $\{\mathbf{P}(X_\delta > x) - \mathbf{P}(X_0 > x)\}/\delta^2 \to \sum_{k=1}^\infty w_k \mu_k^2 h_k(x)$ as $\delta \to 0$, where $h_k$ is a density whose associated characteristic function $\hat{f}(t,0)/(1 - 2iw_k t) = (1 - 2iw_k t)^{-1} \prod_{j=1}^\infty (1 - 2iw_j t)^{-1/2}$ is that of $X_0 + w_k \chi_2^2$, in which the summands are taken to be independent.*

PROOF. It follows from the definition of $X_\delta$ that

$$\hat{f}(t,\delta) = \mathrm{E}(e^{itX_\delta}) = \prod_{k=1}^\infty \hat{\eta}(w_k t, \delta\mu_k) = \xi(t) e^{-2\delta^2 t^2 \kappa_1(t)} e^{i\kappa_2(t) + i\delta^2 \kappa_3(t)},$$

where

$$\xi(t) = \prod_{k=1}^\infty (1 + 4t^2 w_k^2)^{-1/4}, \qquad \kappa_1(t) = \sum_{k=1}^\infty w_k^2 \mu_k^2/(1 + 4t^2 w_k^2),$$

$$\kappa_2(t) = \frac{1}{2}\sum_{k=1}^\infty \arctan(2tw_k), \qquad \kappa_3(t) = t\sum_{k=1}^\infty w_k \mu_k^2/(1 + 4t^2 w_k^2).$$

Note that $\xi(t)$ and $t^2\xi(t)$ are integrable, $\kappa_1$ is bounded, $\kappa_i(t)/t$ is bounded for $i \in \{2,3\}$ and that as $t \to 0$, $\kappa_2(t)/t \to 1/36$ and $\kappa_3(t)/t \to \sum_{k=1}^\infty w_k \mu_k^2$.

Next, from the Gil-Pelaez representation, one has

$$\mathbf{P}(X_\delta > x) - \mathbf{P}(X_0 > x) = \frac{1}{2\pi}\int_{-\infty}^{+\infty} t^{-1}\mathrm{Im}\{e^{-itx}\hat{f}(t,\delta) - e^{-itx}\hat{f}(t,0)\}\,dt.$$

Note that

$$t^{-1}\mathrm{Im}\{e^{-itx}\hat{f}(t,\delta)\} = t^{-1}\xi(t)e^{-2\delta^2 t^2 \kappa_1(t)}\sin\{\kappa_2(t) + \delta^2\kappa_3(t) - tx\}.$$

As a result, $(\delta^2 t)^{-1} \operatorname{Im}\{e^{-itx}\hat{f}(t,\delta) - e^{-itx}\hat{f}(t,0)\}$ can be decomposed as $A_1(t,\delta)t^2\xi(t) + A_2(t,\delta)\xi(t)$, where

$$A_1(t,\delta) = (\delta^2 t^3)^{-1}\{e^{-2\delta^2 t^2 \kappa_1(t)} - 1\}\sin\{\kappa_2(t) + \delta^2\kappa_3(t) - tx\}$$

and

$$A_2(t,\delta) = (\delta^2 t)^{-1}[\sin\{\kappa_2(t) + \delta^2\kappa_3(t) - tx\} - \sin\{\kappa_2(t) - tx\}].$$

Both terms are bounded and converge as $\delta \to 0$. Their limits are $A_1(t,0) = -2t^{-1}\kappa_1(t)\sin\{\kappa_2(t) - tx\}$ and $A_2(t,0) = t^{-1}\kappa_3(t)\cos\{\kappa_2(t) - tx\}$, respectively. An application of Lebesgue's dominated convergence theorem thus yields

$$\lim_{\delta \to 0} \delta^{-2} \int_{-\infty}^{\infty} t^{-1} \operatorname{Im}\{e^{-itx}\hat{f}(t,\delta) - e^{-itx}\hat{f}(t,0)\} dt = \int_{-\infty}^{\infty} \psi(t,x) dt,$$

where $\psi(t,x) = \xi(t)[t^{-1}\kappa_3(t)\cos\{\kappa_2(t) - tx\} - 2t\kappa_1(t)\sin\{\kappa_2(t) - tx\}]$.

It is easy to check that $\psi$ can also be expressed as

$$\psi(t,x) = \sum_{k=1}^{\infty} w_k \mu_k^2 \operatorname{Re}\{e^{-itx}\hat{f}(t,0)(1 - 2itw_k)^{-1}\},$$

where $\operatorname{Re}(z)$ stands for the real part of any complex number $z$.

Since $\xi$ is integrable, so is $\hat{f}(t,0)/(1 - 2itw_k)$ and hence

$$\frac{1}{2\pi} \int_{-\infty}^{\infty} \psi(t,x) dt = \sum_{k=1}^{\infty} \frac{w_k \mu_k^2}{\pi} \int_0^{\infty} \operatorname{Re}\{e^{-itx}\hat{f}(t,0)(1 - 2itw_k)^{-1}\} dt$$

$$= \sum_{k=1}^{\infty} w_k \mu_k^2 h_k(x),$$

where $h_k$ is the density of $X_0 + w_k \chi_2^2$, whose summands are taken to be independent. This completes the proof. $\square$

C. GENEST
DÉPARTEMENT DE MATHÉMATIQUES
ET DE STATISTIQUE
UNIVERSITÉ LAVAL
QUÉBEC, QUÉBEC
CANADA G1K 7P4
E-MAIL: Christian.Genest@mat.ulaval.ca

J.-F. QUESSY
DÉPARTEMENT DE MATHÉMATIQUES
ET D'INFORMATIQUE
UNIVERSITÉ DU QUÉBEC À TROIS-RIVIÈRES
CASIER POSTAL 500
TROIS-RIVIÈRES, QUÉBEC
CANADA G9A 5H7
E-MAIL: Jean-Francois.Quessy@uqtr.ca

B. Rémillard
Service de l'enseignement
  des méthodes quantitatives de gestion
3000, chemin de la Côte-Sainte-Catherine
HEC Montréal
Montréal, Québec
Canada H3T 2A7
E-mail: Bruno.Remillard@hec.ca